# LIE SUPERALGEBRAS GRADED BY THE ROOT SYSTEM A(m,n)

Georgia Benkart[1]
Alberto Elduque[2]

February 25, 2002

Abstract. We determine the Lie superalgebras that are graded by the root systems of the basic classical simple Lie superalgebras of type $\mathrm{A}(m,n)$.

## §1. Introduction

Our investigation is a natural extension of work on the problem of classifying Lie algebras graded by finite root systems. Many important classes of Lie algebras such as the affine and toroidal Lie algebras and various generalizations of them, such as the intersection matrix Lie algebras of Slodowy [S], which arise in the study of singularities, or the extended affine Lie algebras of [AABGP], exhibit a grading by a finite (possibly nonreduced) root system $\Delta$. The formal definition, as first given by Berman and Moody in [BM], depends on a finite-dimensional split simple Lie algebra $\mathfrak{g}$ over a field $\mathbb{F}$ of characteristic zero having a root space decomposition $\mathfrak{g} = \mathfrak{h} \oplus \bigoplus_{\mu\in\Delta} \mathfrak{g}_\mu$ relative to a split Cartan subalgebra $\mathfrak{h}$. Such a Lie algebra $\mathfrak{g}$ is an analogue over $\mathbb{F}$ of a finite-dimensional complex simple Lie algebra.

**Definition 1.1.** *A Lie algebra $L$ over $\mathbb{F}$ is* **graded by the (reduced) root system $\Delta$** *or is* **$\Delta$-graded** *if*

($\Delta$G1) *$L$ contains as a subalgebra a finite-dimensional split simple Lie algebra $\mathfrak{g} = \mathfrak{h} \oplus \bigoplus_{\mu\in\Delta} \mathfrak{g}_\mu$ whose root system is $\Delta$ relative to a split Cartan subalgebra $\mathfrak{h} = \mathfrak{g}_0$;*

($\Delta$G2) *$L = \bigoplus_{\mu\in\Delta\cup\{0\}} L_\mu$, where $L_\mu = \{x \in L \mid [h,x] = \mu(h)x \text{ for all } h \in \mathfrak{h}\}$ for $\mu \in \Delta\cup\{0\}$; and*

($\Delta$G3) $L_0 = \sum_{\mu\in\Delta}[L_\mu, L_{-\mu}]$.

2000 *Mathematics Subject Classification.* Primary 17A70.

[1]Support from National Science Foundation Grants #DMS–9810361 (at the Mathematical Sciences Research Institute, Berkeley) and #DMS–9970119 is gratefully acknowledged.
[2]Supported by the Spanish DGI (BFM2001-3239-C03-03).

The second condition in Definition 1.1 describes the grading, while the first provides a uniformity to the structure and enables the representation theory of $\mathfrak{g}$ to be used to study $L$. Condition ($\Delta$G3) insures that the space $L_0$ is connected to the root spaces. Without it, any central ideal could be added to $L$. If only ($\Delta$G1) and ($\Delta$G2) are assumed, then the subalgebra $L' = \Big(\sum_{\mu\in\Delta}[L_\mu, L_{-\mu}]\Big) \oplus \bigoplus_{\mu\in\Delta} L_\mu$ is $\Delta$-graded.

There is a parallel notion of a Lie algebra graded by the nonreduced root system $\mathrm{BC}_r$ introduced and studied in [ABG2] (see also [BS] for the $\mathrm{BC}_1$-case). The Lie algebras graded by finite root systems (both reduced and nonreduced) decompose relative to the adjoint action of $\mathfrak{g}$ into a direct sum of finite-dimensional irreducible $\mathfrak{g}$-modules. There is one possible isotypic component corresponding to each root length and one corresponding to 0 (the sum of the trivial $\mathfrak{g}$-modules). Thus, for the simply-laced root systems only adjoint modules and trivial modules occur. For the doubly-laced root systems, copies of the module having the highest short root as its highest weight also can occur. For type $\mathrm{BC}_r$, there are up to four isotypic components, except when the grading subalgebra $\mathfrak{g}$ has type $\mathrm{D}_2 \cong \mathrm{A}_1\times \mathrm{A}_1$, where there are five possible isotypic components. The complexity increases with the number of isotypic components. The representation theory of $\mathfrak{g}$ is an essential ingredient in the classification of the Lie algebras graded by finite root systems, which has been accomplished in the papers [BM], [BZ], [N], [ABG1], [ABG2], [BS].

Our focus here and in [BE1], [BE2] is on Lie superalgebras graded by the root systems of the finite-dimensional basic classical simple Lie superalgebras $\mathrm{A}(m,n)$, $\mathrm{B}(m,n)$, $\mathrm{C}(n)$, $\mathrm{D}(m,n)$, $\mathrm{D}(2,1;\alpha)$ ($\alpha \in \mathbb{F}\backslash\{0,-1\}$), $\mathrm{F}(4)$, and $\mathrm{G}(3)$. (A standard reference for results on simple Lie superalgebras is Kac's seminal paper [K1].)

Let $\mathfrak{g}$ be a finite-dimensional split simple basic classical Lie superalgebra over a field $\mathbb{F}$ of characteristic zero with root space decomposition $\mathfrak{g} = \mathfrak{h} \oplus \bigoplus_{\mu\in\Delta} \mathfrak{g}_\mu$ relative to a split Cartan subalgebra $\mathfrak{h}$. Thus, $\mathfrak{g}$ is an analogue over $\mathbb{F}$ of a complex simple Lie superalgebra whose root system $\Delta$ is of type $\mathrm{A}(m,n)$ ($m \geq n \geq 0$, $m+n \geq 1$), $\mathrm{B}(m,n)$ ($m \geq 0$, $n \geq 1$), $\mathrm{C}(n)$ ($n \geq 3$), $\mathrm{D}(m,n)$ ($m \geq 2$, $n \geq 1$), $\mathrm{D}(2,1;\alpha)$ ($\alpha \in \mathbb{F} \setminus \{0,-1\}$), $\mathrm{F}(4)$, and $\mathrm{G}(3)$. These Lie superalgebras can be characterized by the properties of being simple, having reductive even part, and having a nondegenerate even supersymmetric bilinear form. Imitating Definition 1.1, we say

**Definition 1.2.** *([BE1, Defn. 2.1]) A Lie superalgebra $L$ over $\mathbb{F}$ is* **graded by the root system $\Delta$** *or is* **$\Delta$-graded** *if*

(i) *$L$ contains as a subsuperalgebra a finite-dimensional split simple basic classical Lie superalgebra $\mathfrak{g} = \mathfrak{h} \oplus \bigoplus_{\mu\in\Delta} \mathfrak{g}_\mu$ whose root system is $\Delta$ relative to a split Cartan subalgebra $\mathfrak{h} = \mathfrak{g}_0$;*

(ii) *($\Delta$G2) and ($\Delta$G3) of Definition 1.1 hold for $L$ relative to the root system $\Delta$.*

The B$(m,n)$-graded Lie superalgebras were determined in [BE1]. These Lie superalgebras differ from others because of their complicated structure and most closely resemble the Lie algebras graded by the nonreduced root systems $\mathrm{BC}_r$. The $\Delta$-graded Lie superalgebras for $\Delta = \mathrm{C}(n)$, $\mathrm{D}(m,n)$, $\mathrm{D}(2,1;\alpha)\,(\alpha \in \mathbb{F} \setminus \{0,-1\})$, F(4), and G(3) were fully described in [BE2]. Therefore, the only remaining case is that of $\Delta = \mathrm{A}(m,n)$ $(m \geq n \geq 0,\, m+n \geq 1)$. Complete results (Theorem 3.10 and Corollary 3.12 below) will be established here for the case $m \neq n$. The Lie superalgebras graded by the root system A$(n,n)$ are truly exceptional for several reasons, and we determine only those A$(n,n)$-graded Lie superalgebras that are completely reducible $\mathfrak{g}$-modules (in Theorem 4.2 and Corollary 4.4).

Our "modus operandi" for investigating the A$(m,n)$-graded Lie superalgebras will mimic that used in our previous papers, whose main steps can be summarized as follows:

**Procedure 1.3.**

1. The determination of the finite-dimensional irreducible $\mathfrak{g}$-modules whose nonzero weights relative to the Cartan subalgebra $\mathfrak{h}$ are roots.
2. The proof of the complete reducibility of any $\Delta$-graded Lie superalgebra $L$ as a module for $\mathfrak{g}$.
3. The computation of $\mathrm{Hom}_{\mathfrak{g}}(V \otimes W, X)$ for any $\mathfrak{g}$-modules $V, W, X$ in Step 1.
4. The determination of the multiplication in any $\Delta$-graded $L$ by combining the previous information.

For A$(m,n)$-graded Lie superalgebras, Step 1 is achieved by arguments similar to those used in the C$(n)$-case of [BE2]. The complete reducibility in Step 2 can be proved along the same lines as in [BE2] provided $m \neq n$. However, for $m = n$, it is no longer true that an A$(n,n)$-graded Lie superalgebra is completely reducible as a module for its grading subsuperalgebra $\mathfrak{g}$, and this is where many complications arise. This case will be considered in the last section of the paper.

The distinctive feature of the A$(m,n)$-graded Lie superalgebras, $m \neq n$, is that $\mathrm{Hom}_{\mathfrak{g}}(\mathfrak{g} \otimes \mathfrak{g}, \mathfrak{g})$ is two-dimensional here, while in all the other cases it is spanned by the Lie bracket.

Steps 1 and 2 above will give that any A$(m,n)$-graded Lie superalgebra $L$ for $m \neq n$, when viewed as a module for the grading subsuperalgebra $\mathfrak{g}$, is a direct sum of two types of irreducible $\mathfrak{g}$-modules – adjoint and trivial ones. By collecting isomorphic summands, we may assume that there are $\mathbb{F}$-vector superspaces $A$ and $D$ so that

$$L = (\mathfrak{g} \otimes A) \oplus D, \tag{1.4}$$

where $D$ is the sum of all the trivial $\mathfrak{g}$-modules. The multiplication on $L$ makes $A$ a superalgebra, while $D$ is a Lie subsuperalgebra acting as superderivations of $A$.

The problem of classifying the A$(m,n)$-graded Lie superalgebras reduces to one of determining the possibilities for $A$ and $D$ and of finding the multiplication. Step 3 will show that $\mathrm{Hom}_{\mathfrak{g}}(\mathfrak{g}\otimes\mathfrak{g},\mathfrak{g})$ is two-dimensional, while $\mathrm{Hom}\mathfrak{g}(\mathfrak{g}\otimes\mathfrak{g},\mathbb{F})$ is spanned by the supertrace. As a result, the product in $L$ must be that given in (3.1) below. The Jacobi identity will impose several restrictions on $A$ and $D$, which eventually will lead to the classification of such Lie superalgebras $L$ (Step 4).

Throughout the paper, $\mathbb{F}$ will denote a fixed but arbitrary field of characteristic zero. Unadorned tensor products will be assumed to be over $\mathbb{F}$.

## §2. The $\mathfrak{g}$-module structure of A$(m,n)$-graded Lie superalgebras ($m>n$)

The following result, proved in [BE1, Lemma 2.2], plays a key role in examining $\Delta$-graded Lie superalgebras.

**Lemma 2.1.** *Let $L$ be a $\Delta$-graded Lie superalgebra, and let $\mathfrak{g}$ be its grading subsuperalgebra. Then $L$ is locally finite as a module for $\mathfrak{g}$.*

As a consequence, each element of a $\Delta$-graded Lie superalgebra $L$, in particular each weight vector of $L$ relative to the Cartan subalgebra $\mathfrak{h}$ of $\mathfrak{g}$, generates a finite-dimensional $\mathfrak{g}$-module. Such a finite-dimensional module has a $\mathfrak{g}$-composition series whose irreducible factors have weights which are roots of $\mathfrak{g}$ or 0. Next we determine (Step 1 in Procedure 1.3) which finite-dimensional irreducible $\mathfrak{g}$-modules have nonzero weights that are roots of $\mathfrak{g}$.

Throughout we will identify the split simple Lie superalgebra $\mathfrak{g}$ of type A$(m,n)$, $m>n\geq 0$, with the special linear Lie superalgebra $\mathfrak{sl}_{m+1,n+1}$. For simplicity of notation, set $p=m+1$ and $q=n+1$, so that $\mathfrak{g}=\mathfrak{sl}_{p,q}$, $p>q\geq 1$.

The diagonal matrices in $\mathfrak{g}$ form a Cartan subalgebra $\mathfrak{h}$, and the corresponding even and odd roots and a system of simple roots of $\mathfrak{g}$ are given by [K1, Sec. 2]:

$$
\begin{aligned}
\Delta_{\bar 0}&=\{\varepsilon_i-\varepsilon_j \mid 1\leq i<j\leq p\}\cup\{\delta_r-\delta_s\mid 1\leq r<s\leq q\},\\
\Delta_{\bar 1}&=\{\pm(\varepsilon_i-\delta_r)\mid 1\leq i\leq p,\ 1\leq r\leq q\},\\
\Pi&=\{\varepsilon_1-\varepsilon_2,\dots,\varepsilon_{p-1}-\varepsilon_p,\varepsilon_p-\delta_1,\delta_1-\delta_2,\dots,\delta_{q-1}-\delta_q\},
\end{aligned}
\tag{2.2}
$$

where for $h=\mathrm{diag}(a_1,\dots,a_p,b_1,\dots,b_q)\in\mathfrak{h}$, $\varepsilon_i(h)=a_i$ and $\delta_r(h)=b_r$ for any $i,r$. The corresponding Cartan matrix is

$$
\left(\begin{array}{ccccccccc}
 & & & & 0 & & & & \\
 & & & & \vdots & & & & \\
 & A_m & & & & & 0 & & \\
 & & & & 0 & & & & \\
 & & & & -1 & & & & \\
0 & \dots & 0 & -1 & 0 & 1 & 0 & \dots & 0 \\
 & & & & -1 & & & & \\
 & & & & 0 & & & & \\
 & 0 & & & \vdots & & A_n & & \\
 & & & & 0 & & & &
\end{array}\right) \tag{2.3}
$$

(if $n = 0$ it is just the $(m+1)\times(m+1)$ upper left corner above), where

$$
A_r = \begin{pmatrix}
2 & -1 & & & \\
-1 & & & & \\
 & & \ddots & & \\
 & & & -1 & \\
 & & -1 & 2 & -1 \\
 & & & -1 & 2
\end{pmatrix}
$$

is the $r \times r$ Cartan matrix correponding to $\mathfrak{sl}_{r+1}$. In terms of the standard matrix units $\{E_{i,j} \mid 1 \leq i,j \leq p+q\}$, the coroots $h_1, \dots, h_{m+n+1}$ have the following expressions:

$$
\begin{aligned}
&h_i = E_{i,i} - E_{i+1,i+1} \qquad\qquad (1 \leq i \leq m+n+1,\ i \neq m+1 = p)\\
&h_p = E_{p,p} + E_{p+1,p+1}.
\end{aligned}
$$

All this is valid too in case $m = n \geq 1$, although there the $h_i$'s are the classes of the elements above modulo the center (and hence they are linearly dependent).

According to [K1, Prop. 2.3]), $\omega \in \mathfrak{h}^* = \mathrm{Hom}_{\mathbb{F}}(\mathfrak{h}, \mathbb{F})$ is a *dominant weight* (hence the highest weight of some finite-dimensional irreducible $\mathfrak{g}$-module) if and only if $\omega(h_i) \in \mathbb{Z}_{\geq 0}$ for all $i \neq p$, and $\omega(h_p) \in \mathbb{Z}$. Therefore, the roots that are dominant weights are $\varepsilon_1 - \varepsilon_p$, $\delta_1 - \delta_q$ (this one does not appear if $q = 1$, that is, if $n = 0$), $\delta_1 - \varepsilon_p$ and $\varepsilon_1 - \delta_q$ (the highest root of $\mathfrak{g}$). Consequently, these are the candidate highest weights for irreducible $\mathfrak{g}$-modules occurring in A$(m,n)$-graded Lie superalgebras.

Now the Lie superalgebra $\mathfrak{g}$ has a $\mathbb{Z}$-gradation, $\mathfrak{g} = \mathfrak{g}_{-1} \oplus \mathfrak{g}_0 \oplus \mathfrak{g}_1$ with $\mathfrak{g}_{\bar{0}} = \mathfrak{g}_0$ and $\mathfrak{g}_{\bar{1}} = \mathfrak{g}_{-1} \oplus \mathfrak{g}_1$, which corresponds to the partition of the $(p+q)\times(p+q)$-matrices into blocks of size $p \times p$, $p \times q$, $q \times p$ and $q \times q$ (if $p = q$ one has to factor out the one-dimensional center). The Lie algebra $\mathfrak{g}_0$ consists of the block diagonal matrices, and $\mathfrak{g}_1$ (respectively $\mathfrak{g}_{-1}$) consists of the block strictly upper (resp. lower) triangular matrices. Kac [K2, Sec. 2] has shown for a finite-dimensional irreducible $\mathfrak{g}$-module $V = V(\Lambda)$ that $V' = \{x \in V \mid \mathfrak{g}_1.x = 0\}$ is an irreducible $\mathfrak{g}_0$-submodule of highest weight $\Lambda$, and $V$ is a quotient of the induced module $\mathcal{U}(\mathfrak{g}) \otimes_{\mathcal{U}(\mathfrak{g}_0 \oplus \mathfrak{g}_1)} V'$, which as a vector space is isomorphic to $\mathfrak{U}(\mathfrak{g}_{-1}) \otimes V'$ (where $\mathfrak{U}(\ )$ denotes the universal enveloping algebra). Thus, the weights of $V$ are of the form $\omega + \nu$, where $\omega$ is a weight of the $\mathfrak{g}_0$-module $V'$ and $\nu$ is a weight of $\mathfrak{U}(\mathfrak{g}_{-1})$. Hence $\nu$ is either 0 or a sum of roots of the form $\delta_r - \varepsilon_i$.

Assume that $\Lambda$ is the root $\delta_1 - \varepsilon_p$. Then no weight $\omega + \nu$ as above is a root or 0 unless $\nu = 0$. Thus $V = V'$ and $\mathfrak{g}_{-1}.V = 0$. This is a contradiction since $\mathfrak{g}$ is simple and $V$ must be a faithful $\mathfrak{g}$-module. Now assume that $\Lambda$ is either $\varepsilon_1 - \varepsilon_p$ or $\delta_1 - \delta_q$. The same argument gives that the possible $\nu$'s involved are either 0 or roots of the form $\delta_r - \varepsilon_i$, so that $V = V' + \mathfrak{g}_{-1}.V'$. If $V = V'$, we reach a contradiction as before, but if $V \neq V'$, then weights of the form $\delta_r - \varepsilon_i$ appear in $V$. Reasoning with the lowest weight instead of the highest one, we see that necessarily the lowest weight of $V$ has to be $\delta_q - \varepsilon_1$, which corresponds to the adjoint module, whose highest weight is $\varepsilon_1 - \delta_q$, a contradiction. Therefore, the only possibility left is that $\Lambda = \varepsilon_1 - \delta_q$, and $V$ is the adjoint module.

Let us summarize the assertions above in the following:

**Theorem 2.4.** *Let $\mathfrak{g}$ be a split simple Lie superalgebra of type $A(m,n)$, with $m \geq n \geq 0$, $m + n \geq 1$. The only finite-dimensional irreducible $\mathfrak{g}$-modules whose weights relative to the Cartan subalgebra of diagonal matrices (modulo the center if necessary) are either roots or 0 are exactly the adjoint and the trivial modules (possibly with the parity changed).*

For showing complete reducibility, we can adapt the proof of complete reducibility in [BE2, Prop. 3.1] to our setting. But from now on, the condition $m > n$ is necessary.

**Proposition 2.5.** *Let $\mathfrak{g}$ be a split simple Lie superalgebra of type $A(m,n)$, $m > n \geq 0$, with split Cartan subalgebra $\mathfrak{h}$. Assume $V$ is a locally finite $\mathfrak{g}$-module satisfying*

- (i) *$\mathfrak{h}$ acts semisimply on $V$;*
- (ii) *any composition factor of any finite-dimensional submodule of $V$ is isomorphic to the adjoint module $\mathfrak{g}$ or to a trivial module (possibly with the parity changed).*

*Then $V$ is a completely reducible $\mathfrak{g}$-module.*

*Proof.* With the same reductions as in [BE2, Prop. 3.1], it suffices to show that if $X$ is a submodule of $V$ and $Y$ is a submodule of $X$ such that either $Y$ is an adjoint module and $X/Y$ is a trivial one-dimensional module, or if $Y$ is trivial and $X/Y$ is adjoint, then $X \cong Y \oplus X/Y$. Since $m > n$, the superform given by the supertrace: $(x,y) \mapsto \mathfrak{str}(xy)$ is a nondegenerate supersymmetric bilinear form on $\mathfrak{g} = \mathfrak{sl}_{m+1,n+1}$. Let $C$ be the Casimir element corresponding to $\mathfrak{g}$ (see [K1, comments before Prop. 5.2.6]). Then one checks easily that $C$ acts as the nonzero scalar $m - n$ on the adjoint module (and it acts trivially on the trivial module). Hence $X$ is the direct sum of the kernel of the action of $C$ and of its image, and thus it is completely reducible. $\square$

Proposition 2.5 is no longer true if $m = n$, and this is the source of many of the difficulties encountered in that case. To illustrate this, consider the Lie superalgebras $\mathfrak{sl}_{n+1,n+1}$ and $\mathfrak{pgl}_{n+1,n+1}$ as modules over the simple Lie superalgebra $\mathfrak{psl}_{n+1,n+1}$ of type A$(n,n)$. In both instances, there is a composition series with an adjoint and a trivial module, but there is no complete reducibility.

The final preparatory step (number 3 in Procedure 1.3) needed to study A$(m,n)$-graded Lie superalgebras is the computation of the spaces of module homomorphisms $\mathrm{Hom}_{\mathfrak{g}}(\mathfrak{g}\otimes\mathfrak{g},\mathbb{F})$ and $\mathrm{Hom}_{\mathfrak{g}}(\mathfrak{g}\otimes\mathfrak{g},\mathfrak{g})$. It is clear that $\mathrm{Hom}_{\mathfrak{g}}(\mathfrak{g}\otimes\mathfrak{g},\mathbb{F})$ is spanned by the even supersymmetric bilinear form induced by the supertrace: $(x,y)\mapsto \mathfrak{str}(xy)$. As for the latter we have

**Proposition 2.6.** *Let $\mathfrak{g}$ be a split simple classical Lie superalgebra of type $A(m,n)$ ($m>n\geq 0$). Then $\mathrm{Hom}_{\mathfrak{g}}(\mathfrak{g}\otimes\mathfrak{g},\mathfrak{g})$ is two-dimensional and spanned by the Lie bracket and by the map given by $(x,y)\mapsto x*y = xy+yx-\frac{2}{m-n}\mathfrak{str}(xy)I$, for any $x,y\in\mathfrak{g}=\mathfrak{sl}_{m+1,n+1}$.*

*Proof.* The Lie bracket and the symmetrized product $x*y$ are module homomorphisms and linearly independent. Hence it is enough to check that the dimension of $\mathrm{Hom}_{\mathfrak{g}}(\mathfrak{g}\otimes\mathfrak{g},\mathfrak{g})\leq 2$.

Let $\mathfrak{g}=\mathfrak{g}_{-1}\oplus\mathfrak{g}_0\oplus\mathfrak{g}_1$ be the $\mathbb{Z}$-gradation considered above. Here $\mathfrak{g}_0=\mathfrak{g}_{\bar 0}=\mathfrak{sl}_p\oplus\mathfrak{sl}_q\oplus\mathbb{F}c$ ($p=m+1$ and $q=n+1$ as before), where $c$ is the block diagonal matrix $c=\mathrm{diag}\,(\frac{q}{q-p}I_p,\frac{p}{q-p}I_q)$, which acts as the identity on $\mathfrak{g}_1$ and as minus the identity on $\mathfrak{g}_{-1}$. Moreover, $\mathfrak{g}_1$ and $\mathfrak{g}_{-1}$ are contragredient irreducible $\mathfrak{g}_0$-modules. Then the matrix $v=E_{1,p+q}$ is a highest weight vector of $\mathfrak{g}_1$ relative to the Cartan subalgebra $\mathfrak{h}$ of $\mathfrak{g}_0$ with respect to the set of simple (even) roots. Also $w=E_{p+q,1}$ is a lowest weight vector for $\mathfrak{g}_{-1}$. As in [BE2, proof of Lemma 5.1], $v\otimes w$ generates $\mathfrak{g}\otimes\mathfrak{g}$ as a $\mathfrak{g}$-module, so any $\varphi\in\mathrm{Hom}_{\mathfrak{g}}(\mathfrak{g}\otimes\mathfrak{g},\mathfrak{g})$ is determined by $\varphi(v\otimes w)$, which belongs to $\mathfrak{h}$ because $v\otimes w$ has weight 0. In particular, $\varphi$ restricts to a nonzero $\mathfrak{g}_0$-module homomorphism $\mathfrak{g}_1\otimes\mathfrak{g}_{-1}\to\mathfrak{g}_0$. But as a module for $[\mathfrak{g}_0,\mathfrak{g}_0]=\mathfrak{sl}_p\oplus\mathfrak{sl}_q$, $\mathfrak{g}_1=V\otimes W^*$ and $\mathfrak{g}_{-1}=V^*\otimes W$, where $V$ is the natural $p$-dimensional module for $\mathfrak{sl}_p$ and $W$ the natural $q$-dimensional module for $\mathfrak{sl}_q$. Since $\mathrm{Hom}_{\mathfrak{sl}_p}(V\otimes V^*,\mathfrak{sl}_p)$ is one-dimensional (as $V\otimes V^*\cong\mathfrak{gl}_p=\mathbb{F}I_p\oplus\mathfrak{sl}_p$), and the same is true for $q$ if $q\geq 1$, it follows that the dimension of $\mathrm{Hom}_{[\mathfrak{g}_0,\mathfrak{g}_0]}(\mathfrak{g}_1\otimes\mathfrak{g}_{-1},[\mathfrak{g}_0,\mathfrak{g}_0])$ is two if $q\geq 1$ and one if $q=1$ ($n=0$). As a consequence, $\dim\mathrm{Hom}_{\mathfrak{g}}(\mathfrak{g}\otimes\mathfrak{g},\mathfrak{g})\leq 3$ if $q\geq 1$, and it is $\leq 2$ if $q=1$. Thus we are done when $n=0$ or when $n\geq 1$, unless there is a $\varphi\in\mathrm{Hom}_{\mathfrak{g}}(\mathfrak{g}\otimes\mathfrak{g},\mathfrak{g})$ such that $\varphi(\mathfrak{g}_1\otimes\mathfrak{g}_{-1})=\mathbb{F}c$. If such a $\varphi$ exists, then for any $x\in\mathfrak{g}_1$ and $y_{\pm1}\in\mathfrak{g}_{\pm1}$,

$$\varphi(y_1\otimes[x,y_{-1}])=\varphi([x,y_1]\otimes y_{-1})-[x,\varphi(y_1\otimes y_{-1})]\in\mathbb{F}[c,x]=\mathbb{F}x.$$

In particular, for $1 \le i \ne j \le p$, and $r \ge p+1$,

$$\varphi(\mathfrak{g}_1 \otimes E_{i,j}) = \varphi(\mathfrak{g}_1 \otimes [E_{i,r}, E_{r,j}]) \subseteq \mathbb{F}E_{i,r},$$

so $\varphi(\mathfrak{g}_1 \otimes E_{i,j}) = 0$ must hold since $q \ge 2$. In a similar vein, $\varphi(\mathfrak{g}_1 \otimes E_{r,s}) = 0$ for any $p+1 \le r \ne s \le p+q$. As the elements $E_{i,j}$ and $E_{r,s}$ generate $[\mathfrak{g}_0, \mathfrak{g}_0]$ as a $\mathfrak{g}_0$-module, it follows that $\varphi(\mathfrak{g}_1 \otimes [\mathfrak{g}_0, \mathfrak{g}_0]) = 0$. Also, for any $1 \le i \le p$, $\varphi(\mathfrak{g}_1 \otimes [E_{i,p+q}, E_{p+q,i}]) \subseteq \mathbb{F}E_{i,p+q}$. But

$$[E_{2,p+q}, E_{p+q,2}] = E_{2,2} + E_{p+q,p+q} = E_{2,2} - E_{1,1} + E_{1,1} + E_{p+q,p+q},$$

so

$$\begin{aligned}\varphi(\mathfrak{g}_1 \otimes [E_{2,p+q}, E_{p+q,2}]) &\subseteq \varphi(\mathfrak{g}_1 \otimes (E_{2,2} - E_{1,1})) + \varphi(\mathfrak{g}_1 \otimes [E_{1,p+q}, E_{p+q,1}]) \\ &\subseteq 0 + \mathbb{F}E_{1,p+q} = \mathbb{F}E_{1,p+q}\end{aligned}$$

too. Hence $\varphi(\mathfrak{g}_1 \otimes (E_{2,2} + E_{p+q,p+q})) = 0$ and

$$\varphi(\mathfrak{g}_1 \otimes \mathfrak{g}_0) = \varphi(\mathfrak{g}_1 \otimes [\mathfrak{g}_0, \mathfrak{g}_0]) + \varphi(\mathfrak{g}_1 \otimes (E_{2,2} + E_{p+q,p+q})) = 0.$$

In the same way one proves that $\varphi(\mathfrak{g}_0 \otimes \mathfrak{g}_{-1}) = 0$. Because $\mathfrak{g} = \mathfrak{g}_{-1} \oplus \mathfrak{g}_0 \oplus \mathfrak{g}_1$ is the eigenspace decomposition for $\operatorname{ad} c$, it follows that $\varphi(\mathfrak{g}_1 \otimes \mathfrak{g}_1) = 0 = \varphi(\mathfrak{g}_{-1} \otimes \mathfrak{g}_{-1})$ also.

Finally,

$$\varphi(\mathfrak{g}_0 \otimes \mathfrak{g}_1) = \varphi([\mathfrak{g}_{-1}, \mathfrak{g}_1] \otimes \mathfrak{g}_1) \subseteq [\mathfrak{g}_{-1}, \varphi(\mathfrak{g}_1 \otimes \mathfrak{g}_1)] + \varphi(\mathfrak{g}_1 \otimes \mathfrak{g}_0) = 0$$

and also $\varphi(\mathfrak{g}_{-1} \otimes \mathfrak{g}_0) = 0$. Hence

$$\varphi(\mathfrak{g} \otimes \mathfrak{g})_{\bar{1}} = \varphi((\mathfrak{g}_1 + \mathfrak{g}_{-1}) \otimes \mathfrak{g}_0) + \varphi(\mathfrak{g}_0 \otimes (\mathfrak{g}_1 + \mathfrak{g}_{-1})) = 0.$$

But $\varphi(\mathfrak{g} \otimes \mathfrak{g})$ is a nonzero ideal of $\mathfrak{g}$, which is simple. This contradiction shows that no $\varphi$ with $\varphi(\mathfrak{g}_1 \otimes \mathfrak{g}_{-1}) = \mathbb{F}c$ exists to complete the proof. $\square$

## §3. The structure of the A$(m,n)$-graded Lie superalgebras $(m > n)$

The results of Section 2 show that any A$(m,n)$-graded Lie superalgebra $L$ with $m > n \ge 0$ is the direct sum of adjoint and trivial modules (possibly with a change of parity) for the grading subalgebra $\mathfrak{g}$. After collecting isomorphic summands, we may suppose that there are superspaces $A = A_{\bar{0}} \oplus A_{\bar{1}}$ and $D = D_{\bar{0}} \oplus D_{\bar{1}}$ so that $L = (\mathfrak{g} \otimes A) \oplus D$, and a distinguished element $1 \in A_{\bar{0}}$ which allows us to identify the grading subalgebra $\mathfrak{g}$ with $\mathfrak{g} \otimes 1$, Observe first that $D$ is a subsuperalgebra of $L$, since it is the (super)centralizer of $\mathfrak{g}$.

For determining the multiplication on $L$, we may apply the same type of arguments as in [BZ] and [BE1,2]. Indeed, fixing homogeneous basis elements $\{a_i\}_{i\in I}$ of $A$ and choosing $a_i, a_j, a_k$ with $i,j,k \in I$, we see that the projection of the product $[\mathfrak{g}\otimes a_i, \mathfrak{g}\otimes a_j]$ onto $\mathfrak{g}\otimes a_k$ determines an element of $\mathrm{Hom}_{\mathfrak{g}}(\mathfrak{g}\otimes\mathfrak{g},\mathfrak{g})$, which is spanned by the supercommutator and the symmetrized product. Thus, there exist scalars $\xi^k_{i,j}$ and $\vartheta^k_{i,j}$ so that

$$\begin{aligned}
[x\otimes a_i, y\otimes a_j]\Big|_{\mathfrak{g}\otimes A} &= (-1)^{\bar a_i \bar y}\left(\sum_{k\in I}\xi^k_{i,j}[x,y]\otimes a_k + \sum_{k\in I}\vartheta^k_{i,j}x*y\otimes a_k\right)\\
&= (-1)^{\bar a_i \bar y}\left([x,y]\otimes\Big(\sum_{k\in I}\xi^k_{i,j}a_k\Big) + x*y\otimes\Big(\sum_{k\in I}\vartheta^k_{i,j}a_k\Big)\right).
\end{aligned}$$

(Our convention is that $\bar y = e$ whenever $y \in \mathfrak{g}_e$ for $e = \bar 0, \bar 1$, etc.) Defining $\circ : A\times A \to A$ by $a_i \circ a_j = 2\sum_{k\in I}\xi^k_{i,j}a_k$, and $[\,,\,] : A\times A\to A$ by $[a_i,a_j] = 2\sum_{k\in I}\vartheta^k_{i,j}a_k$ and extending them bilinearly, we obtain two products "$\circ$" and "$[\,,\,]$" on $A$. (The factors of 2 are simply for convenience.) Necessarily the first is supercommutative and the second superanticommutative, because the products on $\mathfrak{g}$ and $L$ are superanticommutative. Using similar reasoning and taking into account that $\mathrm{Hom}_{\mathfrak{g}}(\mathfrak{g}\otimes\mathfrak{g},\mathbb{F})$ is spanned by the supertrace, we see that there exist a (super) skew symmetric form $\langle\,|\,\rangle : A\times A\to D$ and an even bilinear map $D\times A\to A$: $(d,a)\mapsto da$ with $d1 = 0$ such that the multiplication in $L$ is given by:

$$\begin{aligned}
&[f\otimes a, g\otimes a'] = (-1)^{\bar a\bar g}\Big([f,g]\otimes\frac12 a\circ a' + f*g\otimes\frac12[a,a'] + \mathfrak{str}(fg)\langle a\mid a'\rangle\Big)\\
&[d, f\otimes a] = (-1)^{\bar d\bar f} f\otimes da,\\
&[d,d'] \quad \text{(the product in } D)
\end{aligned} \tag{3.1}$$

for homogeneous elements $f,g\in\mathfrak{g}$, $a,a'\in A$, $d,d'\in D$. Additionally, $1\circ a = 2a$ for any $a\in A$ and $[1,A] = 0$.

There is a unique unital multiplication $aa'$ on $A$ such that $a\circ a' = aa' + (-1)^{\bar a\bar a'}a'a$ and $[a,a'] = aa' - (-1)^{\bar a\bar a'}a'a$ for any homogeneous elements $a,a'\in A$. When we refer to the algebra $A$, it is this unital multiplication that will be tacitly assumed.

Now the Jacobi superidentity $\mathfrak{J}(z_1,z_2,z_3) = \sum_{\circlearrowleft}(-1)^{\bar z_1\bar z_3}[[z_1,z_2],z_3] = 0$ (cyclic permutation of the homogeneous elements $z_1,z_2,z_3$), when specialized with homogeneous elements $d_1,d_2\in D$ and $f\otimes a\in\mathfrak{g}\otimes A$ shows that $\phi : D\to\mathrm{End}_{\mathbb{F}}(A)$: $\phi(d)(a) = da$, is a representation of the Lie superalgebra $D$. When it is specialized with homogeneous elements $d\in D$ and $f\otimes a$, $g\otimes a'\in\mathfrak{g}\otimes A$, we obtain

$$[d,[f\otimes a, g\otimes a']] = [[d,f\otimes a], g\otimes a'] + (-1)^{\bar d(\bar f+\bar a)}[f\otimes a,[d,g\otimes a']],$$

and using (3.1), we see that this is the same as:

(3.2)
$$
\begin{aligned}
&(-1)^{\bar d(\bar f+\bar g)}(-1)^{\bar a\bar g}\Bigl([f,g]\otimes\frac{1}{2}d(a\circ a')+f*g\otimes\frac{1}{2}d([a,a'])+\mathfrak{str}(fg)[d,\langle a\mid a'\rangle]\Bigr)\\
&\ =(-1)^{\bar d\bar f}(-1)^{(\bar d+\bar a)\bar g}\Bigl([f,g]\otimes\frac{1}{2}(da)\circ a'+f*g\otimes\frac{1}{2}[(da),a']+\mathfrak{str}(fg)\langle da\mid a'\rangle\Bigr)\\
&\ +(-1)^{\bar d(\bar f+\bar a+\bar g)}(-1)^{\bar a\bar g}\Bigl([f,g]\otimes\frac{1}{2}a\circ(da')+f*g\otimes\frac{1}{2}[a,(da')]+\mathfrak{str}(fg)\langle a\mid da'\rangle\Bigr).
\end{aligned}
$$

When $f=E_{1,2}$ and $g=E_{2,1}$, the elements $[f,g]$ and $f*g$ are linearly independent and $\mathfrak{str}(fg)=1$. Hence (3.2) is equivalent to:

(i) $d(a\circ a')=(da)\circ a'+(-1)^{\bar d\bar a}a\circ(da')$,
(ii) $d([a,a'])=[(da),a']+[a,(da')]$,
(iii) $[d,\langle a\mid a'\rangle]=\langle da\mid a'\rangle+(-1)^{\bar d\bar a}\langle a\mid da'\rangle$,

for any homogeneous $d\in D$ and $a,a'\in A$. Items (i) and (ii) can be combined to give

(3.3) $\phi$ is a representation as superderivations: $\phi: D\to \mathrm{Der}_{\mathbb{F}}(A)$,

while (iii) says that

(3.4) $\langle\,|\,\rangle$ is invariant under the action of $D$.

For $z_1\otimes a_1$, $z_2\otimes a_2$, $z_3\otimes a_3\in\mathfrak{g}\otimes A$, the Jacobi superidentity is equivalent to the two relations

(3.5)
$$
\begin{aligned}
0=&\,\mathfrak{str}(z_1z_2z_3)\Bigl(\sum_{\circlearrowleft}(-1)^{\bar a_1\bar a_3}\langle a_1\mid a_2a_3\rangle\Bigr)\\
&-(-1)^{\bar z_2\bar z_3}\mathfrak{str}(z_1z_3z_2)\Bigl(\sum_{\circlearrowleft}(-1)^{(\bar a_1+\bar a_2)\bar a_3}\langle a_1\mid a_3a_2\rangle\Bigr)
\end{aligned}
$$

(3.6)

$$\begin{aligned}
0 = &- \sum_{\circlearrowleft} (-1)^{\bar{z_1}\bar{z_3}+\bar{a_1}\bar{a_3}} z_1z_2z_3 \otimes (a_1,a_2,a_3) \\
&+ \sum_{\circlearrowleft} (-1)^{(\bar{z_1}+\bar{z_2})\bar{z_3}+(\bar{a_1}+\bar{a_2})\bar{a_3}} z_1z_3z_2 \otimes (a_1,a_3,a_2) \\
&- \sum_{\circlearrowleft} (-1)^{\bar{z_1}\bar{z_3}+\bar{a_1}\bar{a_3}} \mathfrak{str}(z_1z_2)z_3 \otimes \Big( \langle a_1 \mid a_2\rangle a_3 - \frac{1}{m-n}[[a_1,a_2],a_3]\Big) \\
&- (-1)^{\bar{z_1}\bar{z_3}} \frac{\mathfrak{str}(z_1z_2z_3)}{m-n} I \otimes \Big( \sum_{\circlearrowleft} (-1)^{\bar{a_1}\bar{a_3}} [a_1,a_2a_3]\Big) \\
&+ (-1)^{(\bar{z_1}+\bar{z_2})\bar{z_3}} \frac{\mathfrak{str}(z_1z_3z_2)}{m-n} I \otimes \Big( \sum_{\circlearrowleft} (-1)^{(\bar{a_1}+\bar{a_2})\bar{a_3}} [a_1,a_3a_2]\Big),
\end{aligned}$$

where $(a_1,a_2,a_3) = (a_1a_2)a_3 - a_1(a_2a_3)$, (the associator). The first corresponds to the $D$-component and the second to the $(\mathfrak{g}\otimes A)$-component. Here (3.6) makes sense inside $\mathfrak{gl}_{p,q}\otimes A \supseteq \mathfrak{sl}_{p,q}\otimes A = \mathfrak{g}\otimes A$, and $I$ is the identity matrix.

Now suppose $z_1 = E_{1,2}$, $z_2 = E_{2,3}$ and $z_3 = E_{3,1}$ in $\mathfrak{g}$. They are all even if $m \geq 2$, while $z_1$ is even and $z_2$ and $z_3$ are odd if $m = 1$ (so $n = 0$). Then $z_iz_{i+1}z_{i+2} = E_{i,i}$, $z_iz_{i+2}z_{i+1} = 0$, $\mathfrak{str}(z_iz_{i+1}) = 0$, $\mathfrak{str}(z_iz_{i+1}z_{i+2}) = \pm 1$ (indices modulo 3). As a result, (3.5) gives:

$$\sum_{\circlearrowleft} (-1)^{\bar{a_1}\bar{a_3}} \langle a_1 \mid a_2a_3\rangle = 0. \tag{3.7}$$

If $m \geq 2$ , $E_{1,1}$, $E_{2,2}$, $E_{3,3}$, and $I$ are linearly independent, so Equation 3.6 gives $(a_1,a_2,a_3) = 0$, that is

$$A \text{ is associative.} \tag{3.8}$$

Then (3.6) becomes simply

$$0 = \sum_{\circlearrowleft} (-1)^{\bar{z_1}\bar{z_3}+\bar{a_1}\bar{a_3}} \mathfrak{str}(z_1z_2)z_3 \otimes \Big( \langle a_1 \mid a_2\rangle a_3 - \frac{1}{m-n}[[a_1,a_2],a_3]\Big),$$

and since $\mathfrak{str}(z_1z_2)z_3$ is not identically 0 on $\mathfrak{g}$,

$$\langle a_1 \mid a_2\rangle a_3 = \frac{1}{m-n}[[a_1,a_2],a_3], \tag{3.9}$$

for any $a_1, a_2, a_3 \in A$.

Now, if $m = 1$, the expression in (3.6) with the $z_i$'s chosen as above is a linear combination of $E_{1,1}$, $E_{2,2}$ and $E_{3,3}$ with coefficients in $A$. The coefficient of $E_{1,1}$ is

$$-(-1)^{\bar{a_1}\bar{a_3}}(a_1, a_2, a_3) - \Big(\sum_{\circlearrowleft}(-1)^{\bar{a_1}\bar{a_3}}[[a_1, a_2], a_3]\Big) = 0,$$

while the coefficient of $E_{3,3}$ is

$$(-1)^{\bar{a_3}\bar{a_2}}(a_3, a_1, a_2) - \Big(\sum_{\circlearrowleft}(-1)^{\bar{a_1}\bar{a_3}}[[a_1, a_2], a_3]\Big) = 0.$$

Hence

$$(a_3, a_2, a_1) = -(-1)^{(\bar{a_1}+\bar{a_2})\bar{a_3}}(a_3, a_1, a_2)$$

for any homogeneous $a_1, a_2, a_3 \in A$. Permuting cyclically twice, we determine that

$$(a_1, a_2, a_3) = -(a_1, a_2, a_3),$$

so that $A$ is associative in this case also, and then (3.9) is satisfied too.

In this way, we have arrived at our main theorem. The last sentence in it is a consequence of condition ($\Delta$G3) in Definition 1.2.

**Theorem 3.10.** *Assume $L = (\mathfrak{g} \otimes A) \oplus D$ is a superalgebra over a field $\mathbb{F}$ of characteristic zero where $\mathfrak{g} = \mathfrak{sl}_{m+1,n+1}$, $m > n \geq 0$, $A$ is unital $\mathbb{F}$-superalgebra, and $D$ is Lie superalgebra, and with multiplication as in (3.1). Then $L$ is a Lie superalgebra if and only if*

- *$A$ is a unital associative superalgebra,*
- *$D$ is a Lie subsuperalgebra of $L$ and $\phi : D \to \mathrm{Der}_{\mathbb{F}}(A)$ ($\phi(d)(a) = da$) is a representation of $D$ as superderivations on the algebra $A$,*
- $[d, \langle a_1 \mid a_2\rangle] = \langle da_1 \mid a_2\rangle + (-1)^{\bar{d}\bar{a}_1}\langle a_1 \mid da_2\rangle$,
- $\sum_{\circlearrowleft}(-1)^{\bar{a_1}\bar{a_3}}\langle a_1 \mid a_2a_3\rangle = 0$,
- $\langle a_1 \mid a_2\rangle a_3 = \frac{1}{m-n}[[a_1, a_2], a_3]$,

*for any homogeneous elements $d \in D$ and $a_1, a_2, a_3 \in A$.*

*Moreover, the $A(m, n)$-graded Lie superalgebras (for $m > n \geq 0$) are exactly these superalgebras with the added constraint that*

$$D = \langle A \mid A\rangle.$$

**Remark 3.11.** Let $A$ be any unital associative superalgebra. Then $\mathrm{ad}_{[A,A]}$ is a subsuperalgebra of $\mathrm{Der}_{\mathbb{F}}(A)$. Consider the Lie superalgebra

$$\mathfrak{L}(A) \overset{\text{def}}{=} (\mathfrak{g} \otimes A) \oplus \mathrm{ad}_{[A,A]},$$

with $\mathfrak{g} = \mathfrak{sl}_{m+1,n+1}$ $(m > n \geq 0)$, with multiplication given by (3.1), with $\mathrm{ad}_{[A,A]}$ in place of $D$, and with $\langle a \mid a' \rangle = \frac{1}{m-n}\mathrm{ad}_{[a,a']}$ for any $a, a' \in A$. Then Theorem 3.10 shows that $\mathfrak{L}(A)$ is an A$(m, n)$-graded Lie superalgebra.

Moreover, for any A$(m, n)$-graded Lie superalgebra $L$ with coordinate superalgebra $A$, Theorem 3.10 implies that $L/Z(L) \cong \mathfrak{L}(A)$. Thus $L$ is a cover of $\mathfrak{L}(A)$ (a central extension of $\mathfrak{L}(A)$ which is perfect, $L = [L, L]$).

Any perfect Lie superalgebra $L$ has a unique (up to isomorphism) universal central extension, which is also perfect, called the *universal covering superalgebra* of $L$. Two perfect Lie superalgebras $L_1$ and $L_2$ are said to be *centrally isogenous* if $L_1/Z(L_1) \cong L_2/Z(L_2)$.

To see how these concepts apply in our case, consider the Lie superalgebra $\mathfrak{gl}_{m+1,n+1}(A)$. This is the Lie superalgebra $\mathfrak{M}_{m+1,n+1}(\mathbb{F}) \otimes A$ of the associative superalgebra of block partitioned matrices of size $\big((m+1)+(n+1)\big) \times \big((m+1)+(n+1)\big)$ tensored (as superalgebras over $\mathbb{F}$) with the associative superalgebra $A$. Its commutator subalgebra $\mathfrak{sl}_{m+1,n+1}(A) = [\mathfrak{gl}_{m+1,n+1}(A), \mathfrak{gl}_{m+1,n+1}(A)]$ is an A$(m, n)$-graded Lie superalgebra (since $A$ is unital) with $A$ as a coordinate superalgebra. Thus, Theorem 3.10 and Remark 3.11 give:

**Corollary 3.12.** *The A$(m, n)$-graded Lie superalgebras with $m > n \geq 0$ are precisely the Lie superalgebras which are centrally isogeneous to the Lie superalgebras $\mathfrak{sl}_{m+1,n+1}(A)$ for $A$ a unital associative superalgebra $A$.*

The universal central extension of the Lie superalgebra $\mathfrak{sl}_{m+1,n+1}(A)$ with $m \neq n$ and $m + n \geq 3$ has been shown to be the Steinberg Lie superalgebra $\mathfrak{st}_{m+1,n+1}(A)$ in [MP].

## §4. A$(n, n)$-graded Lie superalgebras

The situation when $m = n$ is much more involved than the previous one, due to the fact that the complete reducibility result in Proposition 2.5 is no longer valid in this case, as already noted above. However, when $L$ is an A$(n, n)$-graded Lie superalgebra with grading subsuperalgebra $\mathfrak{g}$, in one respect the situation is even simpler than for $m \neq n$, because of the next result:

**Proposition 4.1.** *Let $\mathfrak{g}$ be a split simple classical Lie superalgebra of type $A(n,n)$ $(n>0)$. Then $\mathrm{Hom}_{\mathfrak{g}}(\mathfrak{g}\otimes\mathfrak{g},\mathfrak{g})$ is spanned by the Lie bracket.*

*Proof.* Write $\bar{x}$ for the class of a matrix $x \in \mathfrak{sl}_{p,p}$ modulo the center $FI_{2p}$ $(p = n+1)$. As in the proof of Proposition 2.6, let $v = \bar{E}_{1,2p}$ and $w = \bar{E}_{2p,1}$. Then any $\varphi \in \mathrm{Hom}_{\mathfrak{g}}(\mathfrak{g}\otimes\mathfrak{g},\mathfrak{g})$ is determined by $\varphi(v\otimes w)$, which belongs to the Cartan subalgebra $\mathfrak{h}$, and it is annihilated by $\bar{E}_{i,i+1}$ for $2 \le i \le p-1$ and $p+1 \le i \le 2p-2$. Therefore, $\varphi(v\otimes w)$ is a linear combination of $(p-1)\bar{E}_{1,1} - (\bar{E}_{2,2} + \cdots + \bar{E}_{p,p})$ and $(p-1)\bar{E}_{2p,2p} - (\bar{E}_{p+1,p+1} + \cdots + \bar{E}_{2p-1,2p-1})$, so that $\dim \mathrm{Hom}_{\mathfrak{g}}(\mathfrak{g}\otimes\mathfrak{g},\mathfrak{g}) \le 2$. If it were 2, there would exist a $\varphi \in \mathrm{Hom}_{\mathfrak{g}}(\mathfrak{g}\otimes\mathfrak{g},\mathfrak{g})$ such that $\varphi(v\otimes w) = (p-1)\bar{E}_{1,1} - (\bar{E}_{2,2} + \cdots + \bar{E}_{p,p})$, but since $[E_{p,p+1}, E_{1,2p}] = 0 = [E_{p,p+1}, E_{2p,1}]$, it follows that $0 = [\bar{E}_{p,p+1}, \varphi(v\otimes w)] = [\bar{E}_{p,p+1}, (p-1)\bar{E}_{1,1} - (\bar{E}_{2,2} + \cdots + \bar{E}_{p,p})] = \bar{E}_{p,p+1}$, a contradiction. Consequently, $\dim \mathrm{Hom}_{\mathfrak{g}}(\mathfrak{g}\otimes\mathfrak{g},\mathfrak{g}) = 1$, and hence it is spanned by the Lie bracket. □

Therefore, any A$(n,n)$-graded Lie superalgebra with grading subsuperalgebra $\mathfrak{g}$ acting completely reducibly on $L$ satisfies the hypotheses of [BE2, Lemma 4.1], where $\kappa$ is the supersymmetric form on $\mathfrak{g} = \mathfrak{psl}_{p,p}$ induced by the supertrace on $\mathfrak{sl}_{p,p}$ (we will denote it by $\mathfrak{str}$ too), and as an immediate consequence we obtain the following:

**Theorem 4.2.** *Let $L$ be a $\Delta$-graded Lie superalgebra over $\mathbb{F}$ with grading subsuperalgebra $\mathfrak{g}$ of type $A(n,n)$ $(n>0)$, and assume that $L$ is a completely reducible module over $\mathfrak{g}$. Then there is a unital (super)commutative associative superalgebra $A$ and a superspace $D$ such that $L = (\mathfrak{g}\otimes A)\oplus D$, with multiplication given by*

$$\begin{aligned} &[f\otimes a, g\otimes a'] = (-1)^{\bar{a}\bar{g}}\Big([f,g]\otimes aa' + \mathfrak{str}(fg)\langle a \mid a'\rangle\Big) \\ &[d, L] = 0 \end{aligned} \tag{4.3}$$

*for homogeneous elements $f,g\in\mathfrak{g}$, $a,a'\in A$, $d,d'\in D$; where $\langle\,|\,\rangle : A\times A\to D$ is a super skew symmetric bilinear even form with $\langle A \mid A\rangle = D$ and satisfying the two-cocycle condition, $\sum_{\circlearrowleft}(-1)^{\bar{a}_1\bar{a}_3}\langle a_1a_2, a_3\rangle = 0$.*

**Corollary 4.4.** *A $\Delta$-graded Lie superalgebra $L$ with grading subalgebra $\mathfrak{g}$ corresponding to a root system $\Delta$ of type $A(n,n)$ and acting completely reducibly on $L$ is a covering of a Lie superalgebra $\mathfrak{g}\otimes A$, where $A$ is a unital supercommutative associative superalgebra.*

A precise description of the structure of an arbitrary A$(n,n)$-graded Lie superalgebra is yet to be fully resolved. We can prove that, even though complete reducibility fails, as a module over the grading superalgebra $\mathfrak{g} = \mathfrak{psl}_{n+1,n+1}$, such a Lie superalgebra must be a direct sum of copies of $\mathfrak{gl}_{n+1,n+1}$, $\mathfrak{sl}_{n+1,n+1}$, $\mathfrak{pgl}_{n+1,n+1}$, $\mathfrak{psl}_{n+1,n+1}$, and trivial modules. However, the large number of possibilities for elements in the various spaces $\mathrm{Hom}_{\mathfrak{g}}(V \otimes W, X)$, where $V, W, X$ are among the modules above, makes accomplishing Step 4 in Procedure 1.3 a daunting task in this case.

Department of Mathematics, University of Wisconsin, Madison, Wisconsin 53706
*E-mail address*: benkart@math.wisc.edu

Departamento de Matemáticas, Universidad de Zaragoza, 50009 Zaragoza, Spain
*E-mail address*: elduque@posta.unizar.es